\documentclass[12pt]{amsart}

\usepackage{fullpage}

\usepackage{amsmath,amsfonts,amssymb}

\usepackage{amsthm}

\usepackage{graphicx,color}
\usepackage{pstricks,pst-plot,pst-node}

\usepackage{hyperref}

\usepackage[notref,notcite,final]{showkeys}

\numberwithin{equation}{section}

\newtheorem{theorem}{Theorem}[section]

\newtheorem{lemma}[theorem]{Lemma}

\theoremstyle{definition}
\newtheorem{definition}[theorem]{Definition}
\newtheorem{example}[theorem]{Example}

\theoremstyle{remark}
\newtheorem{remark}[theorem]{Remark}

\def\Le{\reflectbox{L}}

\begin{document}


\author{Alexander Burstein}
\address{Department of Mathematics, Iowa State University, Ames, IA 
50011-2064, USA}
\email{burstein@math.iastate.edu}
\urladdr{http://www.math.iastate.edu/burstein/}
\thanks{Research supported in part by the NSA Young Investigator 
Grant H98230-06-1-0037.}

\date{December 31, 2006}

\title{On some properties of permutation tableaux}

\begin{abstract}
We consider the relation between various permutation statistics and
properties of permutation tableaux. We answer some of the questions
of Steingr\'{\i}msson and Williams~\cite{StWi}, in particular, on
the distribution of the bistatistic of numbers of rows and essential
ones in permutation tableaux. We also consider and enumerate sets of
permutation tableaux related to some pattern restrictions on
permutations.
\end{abstract}

\maketitle

\section{Introduction}

Permutation tableaux are a particular class of \Le-diagrams that
were studied by Postnikov~\cite{Postnikov} and enumerated by
Williams~\cite{Williams}. They are defined as follows. Given a
partition of $\lambda=(\lambda_1\ge\lambda_2\ge\dots\ge\lambda_k>0)$
of an integer $m=\sum_i{\lambda_i}$, its \emph{Young diagram}
$Y_\lambda$ of shape $\lambda$ is a left-justified diagram of $m$
boxes, with $\lambda_i$ boxes in row $i$.

A \emph{permutation tableau} $\mathcal{T}^k_n$, which we will also
call a \emph{1-hinge tableau}, is a partition $\lambda$ whose Young
diagram $Y_\lambda$ is contained in a $k\times(n-k)$ rectangle
aligned with its top and left edges, together with a filling of the
cells of $Y_\lambda$ with 0s and 1s that satisfies the following
properties:
\begin{description}
\item[(column)]
Each column of the rectangle contains at least one $1$.

\item[(1-hinge)]
A cell in $Y_\lambda$ with a $1$ above it in the same column and a
$1$ left of it in the same row must contain a $1$.
\end{description}

The filling satisfying the column and hinge properties is called
\emph{valid}. The 1s in $\mathcal{T}^k_n$ that are topmost in their
columns or leftmost in their rows are called \emph{essential}, and
the remaining 1s, i.e. those that are forced by the 1-hinge
property, are called \emph{induced}. Note that the column property
implies that $Y_\lambda$ must have exactly $n-k$ columns and at
least $k$ rows. Note also that some rows of $\mathcal{T}^k_n$ (as
opposed to columns) may contain all zeros. Removing the hinge
requirement yields the definition of \Le-diagram~\cite{Postnikov}.

Alternatively, a permutation tableau $\mathcal{T}^k_n$ may be
thought of as a filling of $k\times(n-k)$ rectangle with 0s, 1s and
2s such that the cells inside $Y_\lambda$ are filled with 0s and 1s
so as to satisfy the column and 1-hinge properties, and the cells
outside $Y_\lambda$ are filled with 2s.

Properties of permutation tableaux were studied by
Steingr\'{\i}msson and Williams~\cite{StWi}. They gave a simpler
description of a map described by Postnikov~\cite{Postnikov2} that
takes permutation tableaux contained in a $k\times(n-k)$ rectangle
to permutations in $\mathcal{S}_n$ with $k$ weak excedances, and
proved that this map $\Phi$ is a bijection that also preserves many
other statistics. The bijection $\Phi$ will be described in the next
section.

One of the conjectures made in~\cite[Section 7]{StWi} is that the
distribution of permutation tableaux according to the number of
essential 1s is equal to that for number of cycles in permutations,
i.e. it is given by the signless Stirling numbers of the first kind.
Moreover, \cite{StWi} conjectured that the joint distribution of
tableaux according to the number of rows and the number of essential
1s equals that of permutations according to the number of weak
excedances and the number of cycles of a permutation, when written
in standard cycle form (which is the same distribution as that of
permutations according to descents and left-to-right-minima). In
this paper we will give a simple natural bijection on \Le-tableaux
that induces the conjectured bijection above and preserves several
other statistics.

Another bijection $\Psi$ on $\mathfrak{S}_n$ defined in~\cite{StWi}
translates certain statistics on permutations corresponding to
entries of their permutation tableaux (determined by $\Phi$) into
certain linear combinations of generalized permutation patterns. In
particular, $\Psi$ maps permutations avoiding generalized pattern
$2\text{-}31$ to permutation tableaux with exactly one 1 in each
column. Since the number of permutations in $\mathfrak{S}_n$
avoiding $2\text{-}31$ is $C_n$, the $n$th Catalan number, we give a
simple description of permutations whose tableau has a single 1 in
each column in terms of noncrossing partitions.

Finally, we describe the properties of the tableaux of permutations
restricted by some 3-letter patterns and enumerate those of them
that have the maximal number of essential 1s.

\section{Permutations and permutation tableaux}

Here we briefly describe a bijection $\Phi$ from \cite{StWi}. We
also give a different proof that $\Phi$ is a bijection, essentially
recovering the tableaux from the corresponding permutation by
reconstructing the columns of the tableau from left to right (or,
similarly, rows from top to bottom), as opposed to right-to-left
column construction in \cite{StWi} (or a similar bottom-to-top row
construction).

Given a permutation tableau $\mathcal{T}^k_n$, its \emph{diagram}
$D(\mathcal{T}^k_n)$ is defined as follows. The southeast border of
partition $\lambda$ gives a path $P=(P_i)_{i=1}^n$ of length $n$
from the northeast corner of the $k\times(n-k)$ rectangle containing
$\lambda$ to the southwestern corner of that rectangle. Label each
step $P_i$ in $P$ with $i$, for $i\in[n]$ (where
$[n]=\{1,2,\dots,n\}$). Now, given an edge $P_i$, also label with
$i$ the edge $Q_i$ on the opposite end of the row (if $P_i$ is
vertical) or column (if $P_i$ is horizontal) containing the edge
$P_i$. Replace each 1 in $\mathcal{T}^k_n$ with a vertex and delete
all 0s. From each vertex draw edges east and south either to the
closest vertices in the same row and the same column or to the
labels $i$ of some edge $P_i$ in $P$. The resulting picture is the
diagram $D(\mathcal{T}^k_n)$. It is also convenient to consider
$D(\mathcal{T}^k_n)$ together with the edges from the labels $i$ of
edges $Q_i$ on the northwestern boundary of $\mathcal{T}^k_n$ to the
closest (leftmost) vertex in the same row (if $Q_i$ is vertical) or
the closest (rightmost) vertex in the same column (if $Q_i$ is
horizontal). We will denote the resulting diagram
$D'(\mathcal{T}^k_n)$ and call it the \emph{expanded diagram} of
$\mathcal{T}^k_n$. It is also convenient to think of edge labels
$i\in[n]$ as labeling a row or a column between $P_i$ and $Q_i$, and
label each cell in the tableau by the ordered pair of its row and
column labels.

Given a tableau $\mathcal{T}^k_n$ as above, the permutation $\pi =
\Phi(\mathcal{T}^k_n)$ as defined as follows. For each $i \in \{1,
\dots , n\}$, consider a zigzag path in $D'(\mathcal{T}^k_n)$ that
starts at $Q_i$ (going south or east depending on whether $Q_i$ is
horizontal or vertical) and switchings direction between south and
east at every vertex it encounters. If that path terminates at
$P_j$, then we set $\pi(i)=j$. Alternatively, in
$D(\mathcal{T}^k_n)$ we replace the first edge of the path starting
from $Q_i$ by an edge north or west from $P_i$ to the vertex in row
or column $i$ that is farthest from $P_i$ (and closest to $Q_i$).

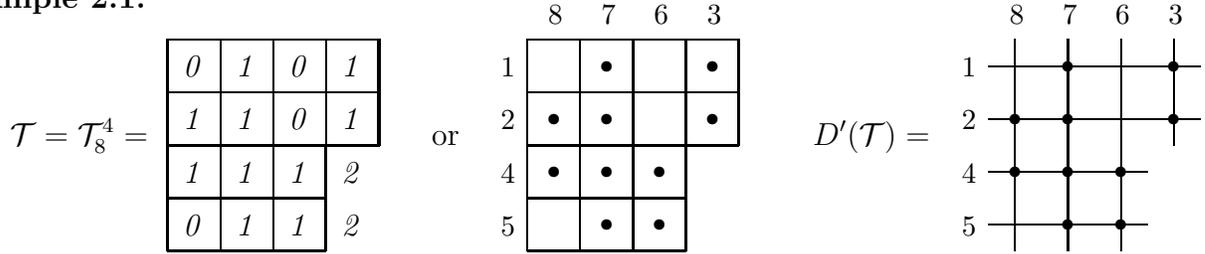
\begin{example}
~\\[-24pt]
\[
\mathcal{T}=\mathcal{T}_{8}^{4}=\mbox{
\begin{picture}(60,60)(20,60)

\setlength{\unitlength}{2pt}

\thinlines
\multiput(10,50)(10,0){4}{\line(0,-1){40}}
\put(50,50){\line(0,-1){20}}
\multiput(10,50)(0,-10){3}{\line(1,0){40}}
\multiput(10,20)(0,-10){2}{\line(1,0){30}}

\thicklines
\put(50,50){\line(0,-1){20}}
\put(50,30){\line(-1,0){10}}
\put(40,30){\line(0,-1){20}}
\put(40,10){\line(-1,0){30}}

\multiput(13,33)(0,-10){2}{\emph{1}}
\multiput(23,43)(0,-10){4}{\emph{1}}
\multiput(33,23)(0,-10){2}{\emph{1}}
\multiput(43,43)(0,-10){2}{\emph{1}}

\multiput(13,43)(0,-30){2}{\emph{0}}
\multiput(33,43)(0,-10){2}{\emph{0}}
\multiput(43,23)(0,-10){2}{\emph{2}}



\end{picture}
}
\hskip 0.5in
\text{or} \quad
\mbox{
\begin{picture}(60,60)(10,60)

\setlength{\unitlength}{2pt}

\thinlines
\multiput(10,50)(10,0){4}{\line(0,-1){40}}
\put(50,50){\line(0,-1){20}}
\multiput(10,50)(0,-10){3}{\line(1,0){40}}
\multiput(10,20)(0,-10){2}{\line(1,0){30}}

\thicklines
\put(50,50){\line(0,-1){20}}
\put(50,30){\line(-1,0){10}}
\put(40,30){\line(0,-1){20}}
\put(40,10){\line(-1,0){30}}

\multiput(15,35)(0,-10){2}{\circle*{2}}
\multiput(25,45)(0,-10){4}{\circle*{2}}
\multiput(35,25)(0,-10){2}{\circle*{2}}
\multiput(45,45)(0,-10){2}{\circle*{2}}

\put(5,43){\small 1}
\put(5,33){\small 2}
\put(5,23){\small 4}
\put(5,13){\small 5}

\put(14,53){\small 8}
\put(24,53){\small 7}
\put(34,53){\small 6}
\put(44,53){\small 3}

\end{picture}
}
\hskip 0.75in
D'(\mathcal{T})=
\mbox{
\begin{picture}(60,60)(5,60)

\setlength{\unitlength}{2pt}



\multiput(15,35)(0,-10){2}{\circle*{2}}
\multiput(25,45)(0,-10){4}{\circle*{2}}
\multiput(35,25)(0,-10){2}{\circle*{2}}
\multiput(45,45)(0,-10){2}{\circle*{2}}

\put(5,43){\small 1}
\put(5,33){\small 2}
\put(5,23){\small 4}
\put(5,13){\small 5}

\put(14,53){\small 8}
\put(24,53){\small 7}
\put(34,53){\small 6}
\put(44,53){\small 3}

\thinlines
\multiput(15,50)(10,0){3}{\line(0,-1){40}}
\put(45,50){\line(0,-1){20}}
\multiput(10,45)(0,-10){2}{\line(1,0){40}}
\multiput(10,25)(0,-10){2}{\line(1,0){30}}

\end{picture}
}
\]
\vskip 0.55in
Following the southeast paths in $D'(\mathcal{T})$ starting with 
each letter from 1 to $n=8$ and switching direction at each dot in 
$D'(\mathcal{T})$, we see that tableau $\mathcal{T}$ corresponds to 
permutation $\pi=\Phi(\mathcal{T})=
\begin{pmatrix}
1&2&3&4&5&6&7&8\\
3&6&1&8&7&4&2&5
\end{pmatrix}
\in\mathfrak{S}_8$. Also, $k=4$ is the number of rows of $\mathcal
{T}$.
\end{example}

It is shown in \cite{StWi} that $\Phi$ is a bijection, $\pi(i)\ge i$
if $P_i$ is vertical and $\pi(i)<i$ if $P_i$ is horizontal, hence
$\pi = \Phi(\mathcal{T}^k_n)$ has $k$ \emph{weak excedances}
(positions $i$ such that $\pi(i)\ge i$) and $n-k$
\emph{deficiencies} (positions $i$ such that $\pi(i)<i$). In fact,
$i$ is a fixed point of $\pi$ (i.e. $\pi(i)=i$) if and only if the
row of $\mathcal{T}^k_n$ labeled $i$ does not contain a 1. We will
refer to $\Phi$ as the standard bijection from permutations tableaux
to permutation and refer to $\Phi(\mathcal{T})$ as the permutation
of $\mathcal{T}$ and to $\Phi^{-1}(\pi)$ as the tableau of $\pi$.

In other words, $\pi$ is a derangement if and only if every row of
$\mathcal{T}=\Phi^{-1}(\pi)$ also must contain a 1. Therefore, we 
may reflect the tableau $\mathcal{T}$ and obtain another
tableau $\mathrm{refl}(\mathcal{T})$ with the same properties.
Moreover, each southeast path $Q_i(\mathcal{T})\to
P_{\pi(i)}(\mathcal{T})$ is thus reflected onto a southeast path
$Q_{n+1-i}(\mathrm{refl}(\mathcal{T}))\to
P_{n+1-\pi(i)}(\mathrm{refl}(\mathcal{T}))$. Thus, if
$\pi=\Phi(\mathcal{T})$ and
$\sigma=\Phi(\mathrm{refl}(\mathcal{T}))$, then
$\sigma(n+1-i)=n+1-\pi(i)$ for all $i\in[n]$, i.e. $\sigma$ is the
reversal of the complement of $\pi$. This answers Open Problem 5 of
\cite{StWi}.

\begin{remark} \label{rem:path-intersect}
Note that each cell with a 0 or 1 in $\mathcal{T}$ is a common point 
of two southeast paths in $D'(\mathcal{T})$ that are uniquely 
determined by the cell, and these paths cross if the cell contains a 
0 and touch, but do not cross, if the cell contains a 1.
\end{remark}

We will need the following result in later sections.

\begin{lemma} \label{lem:path-cross}
For any tableau $\mathcal{T}$, any two southeast paths in 
$D'(\mathcal{T})$ cross at most once, and that intersection must be in the first leg of at least one of the paths (i.e. before reaching the 
first 1). Thus, each 0 in $\mathcal{T}$ corresponds to a unique pair 
of crossing paths in $D'(\mathcal{T})$.
\end{lemma}

\begin{proof}
Suppose that two paths have a common cell $c$ that is in neither of 
their respective first legs. Then one of the paths approaches it 
from a 1 higher in the same column, and the other path approaches it 
from a 1 to the left in the same row. Hence, by the 1-hinge rule the 
cell $c$ also must contain a 1, so these paths can only touch, but 
not intersect, at $c$. 
\end{proof}

Note that each of the two intersecting paths may correspond to an excedance (\textsc{e}) or a non-excedance (\textsc{n}).  
Suppose that two southeast paths cross at a cell $c$ labelled $(i,j)$ (i.e. in row labelled $i$ and column labelled $j$). Then $i<j$. Let $p_{row}$ and $p_{col}$ be the paths entering $c$ from the left and from above, respectively.
 
If $p_{row}$ starts with a south edge $p_{col}$ start with an east edge, then $c$ does not lie on either of their first leg, so these paths cannot cross. Thus, either $p_{row}$ starts with an east edge or $p_{col}$ starts with a south edge, so we may have three types of southeast path intersections in $D'(\mathcal{T})$: \textsc{ee, nn, en}. Let $\mathrm{\textsc{ee}}(\mathcal{T})$, $\textsc{nn}(\mathcal{T})$, $\textsc{en}(\mathcal{T})$ be the sets of corresponding southeast path intersections in $\mathcal{T}$.

If we have an \textsc{en} intersection, then the $p_{col}$ starts south at column labelled $j_{0}>j$ and $p_{row}$ starts east at row labelled $i_{0}$. Then $i_{0}\le i<j\le j_{0}$ (and either $i=i_{0}$ or $j=j_{0}$), so $i_{0}<j_{0}$.

Cells $(i_{0},j_{0})$ with $i_{0}>j_{0}$ that are filled with 2s correspond to \textsc{ne} pairs of southeast paths $p_{row}$ starting at row labelled $i_{0}$ and $p_{col}$ starting at column labelled $j_{0}$, which can never touch or cross (see \cite{StWi}). We denote the set of these cells by $\textsc{ne}(\mathcal{T})$.

\begin{example} \label{ex:ee-nn-en-ne}
~\\[-18pt]
\begin{center}
\begin{picture}(60,60)(0,42)

\setlength{\unitlength}{2pt}

\thinlines
\multiput(10,50)(10,0){4}{\line(0,-1){40}}
\put(50,50){\line(0,-1){20}}
\multiput(10,50)(0,-10){3}{\line(1,0){40}}
\multiput(10,20)(0,-10){2}{\line(1,0){30}}

\thicklines
\put(50,50){\line(0,-1){20}}
\put(50,30){\line(-1,0){10}}
\put(40,30){\line(0,-1){20}}
\put(40,10){\line(-1,0){30}}

\multiput(15,35)(0,-10){2}{\circle*{2}}
\multiput(25,45)(0,-10){4}{\circle*{2}}
\multiput(35,25)(0,-10){2}{\circle*{2}}
\multiput(45,45)(0,-10){2}{\circle*{2}}

\put(5,43){\small 1}
\put(5,33){\small 2}
\put(5,23){\small 4}
\put(5,13){\small 5}

\put(14,53){\small 8}
\put(24,53){\small 7}
\put(34,53){\small 6}
\put(44,53){\small 3}

\put(12,44){\textsc{en}}
\put(12,14){\textsc{ee}}
\put(32,44){\textsc{nn}}
\put(32,34){\textsc{en}}
\put(42,24){\textsc{ne}}
\put(42,14){\textsc{ne}}

\end{picture}
\end{center}
\end{example}

\vskip 0.3in

Several permutation statistics related to corresponding tableaux were also introduced in \cite{StWi}. These statistics essentially count all possible types of pairs of southeast paths ($i\mapsto\pi(i)$, $j\mapsto\pi(j)$) in the tableau corresponding to a given permutation. 
\begin{equation} \label{eq:def-alignments-crossings}
\begin{split}
 A_{EE}(\pi)&=|\{(i,j)\ |\ j<i\le\pi(i)<\pi(j)\}|\\
 A_{NN}(\pi)&=|\{(i,j)\ |\ \pi(j)<\pi(i)<i<j\}|\\
 A_{EN}(\pi)&=|\{(i,j)\ |\ j\le\pi(j)<\pi(i)<i\}|\\
 A_{NE}(\pi)&=|\{(i,j)\ |\ \pi(i)<i<j\le\pi(j)\}|\\
 C_{EE}(\pi)&=|\{(i,j)\ |\ j<i\le\pi(j)<\pi(i)\}|\\
 C_{NN}(\pi)&=|\{(i,j)\ |\ \pi(i)<\pi(j)<i<j\}|
\end{split}
\end{equation}
Then it is easy to see that
\begin{equation} \label{eq:alignments-zeros}
\begin{split}
A_{EE}(\pi)&=|\textsc{ee}(\mathcal{T})|\\
A_{NN}(\pi)&=|\textsc{nn}(\mathcal{T})|\\
A_{EN}(\pi)&=|\textsc{en}(\mathcal{T})|\\
A_{NE}(\pi)&=|\textsc{ne}(\mathcal{T})|=\#\text{2s}(\mathcal{T})
\end{split}
\end{equation}
so the discussion above gives a direct bijective proof that
\begin{equation} \label{eq:count-zeros}
A_{EE}(\pi)+A_{NN}(\pi)+A_{EN}(\pi)=\#\text{0s}(\mathcal{T}),
\end{equation}
which was shown in \cite{StWi} by a more complicated argument. We also note that \cite{StWi} shows that
\[
C_{EE}(\pi)+C_{NN}(\pi)=\#\text{nontop 1s}(\mathcal{T})
\]

Now we will introduce a bit of terminology.
\begin{definition} \label{def:wex}
If $\pi(i)\ge i$ is a weak excedance (\emph{wex}), we will call $\pi
(i)$ a \emph{weak excedance top (wex top)} of $\pi$, and call $i$ a
\emph{weak excedance bottom (wex bottom)} of $\pi$. If $\pi(i)<i$,
we will call $\pi(i)$ a \emph{non-weak-excedance bottom (nonwex
bottom)} of $\pi$, and call $i$ a \emph{non-weak-excedance top
(nonwex top)} of $\pi$.
\end{definition}

We will now describe a way to recover $\mathcal{T}^k_n$ from
$\pi=\Phi(\mathcal{T}^k_n)\in \mathfrak{S}_n$ starting from the
leftmost column. Let $m$ be the largest position of a non-fixed
point in $\pi$. Then $\pi(m)<m$, and $\pi(i)=i$ for all $i>m$. Since
it is trivial to recover the edges $P_i$ for $i>m$, we may assume
that without loss of generality that $m=n$ (so $P_n$ is a horizontal
edge).

\begin{theorem} \label{thm:left-column-dots}
Assume that $\pi\in \mathfrak{S}_n$ is such that $\pi(n)<n$, and let
$\mathcal{T}=\Phi^{-1}(\pi)$. Suppose that the dots in the leftmost
column of $\mathcal{T}$ are in rows labeled $i_1<i_2<\dots<i_r$.
Furthermore, let $\mathcal{T}'$ be the tableau obtained by removing
the leftmost column of $\mathcal{T}$ and replacing $P_n$ with a
vertical edge (so that $n$ becomes a fixed point), and let
$\pi'=\Phi(\mathcal{T}')$. Then
$\pi=\pi'\circ(i_1\,i_2\,\dots\,i_r\,n)$,
$\pi(n)<\pi(i_1)<\dots\pi(i_r)$ and $\pi(i_1)<\dots\pi(i_r)$ are the
successive non-fixed-point left-to-right maxima of the subsequence
of $\pi$ consisting of values greater than $\pi(n)$.
\end{theorem}

\begin{proof}
For $1\le j<r$, the first three steps of the southeast path from
$Q_{i_j}$ to $P_{\pi(i_j)}$ are east from $Q_{i_j}$ to cell
$(i_j,n)$, then south to $(i_{j+1},n)$, then east from
$(i_{j+1},n)$. Thus, $\pi(i_j)=\pi'(i_{j+1})$. Similarly, the path
starting from $Q_{i_r}$ goes east to $(i_r,n)$, then south to $P_n$,
so $\pi(i_r)=n=\pi'(n)$. Likewise, the path from $Q_n$ starts south
to $(i_1,n)$, then turns east, so $\pi(n)=\pi'(i_1)$. Thus,
$\pi=\pi'\circ(i_1\,i_2\,\dots\,i_r\,n)$ as claimed.

Note that paths $p(i_j):i_j\to\pi(i_j)$ and
$p(i_{j+1}):i_{j+1}\to\pi(i_{j+1})$ meet at a vertex in cell
$(i_{j+1},n)$, where $p(i_j)$ enters it traveling south and leaves
east while $p(i_{j+1})$ enters it traveling east and leaves south.
Hence, we can see by induction that at each row $p(i_{j+1})$ is to
the east of $p(i_j)$ and at each column $p(i_j)$ is to the south of
$p(i_{j+1})$, so if $p(i_j)$ and $p(i_{j+1})$ meet again at a cell
it must contain an induced 1, so they cannot cross. It follows from
that $\pi(i_j)<\pi(i_{j+1})$ for all $j<r$, and similarly that
$\pi(n)<\pi(i_1)$.

Now let $l\in[n]$ be a row label such that $i_{j-1}<l<i_{j}$ (if
$j=1$, we simply let $l<i_1$). Then the path from $Q_l$ starts east
and either continues to $P_l$ or first turns south at a cell $(l,m)$
for some $m<n$. Since $(l,m)$ (if it exists) is northeast of
$(i_{j+1},n)$, the same argument as before applies to show that
$\pi(i_j)>\pi(l)$. Therefore, $\pi(i_{j+1})$ is the leftmost value
to the right of $\pi(i_j)$ that is greater than $\pi(i_j)$, so the
theorem follows.
\end{proof}

Note that largest entry $n$ in the (increasing) cycle above gives
the label of the leftmost column, while the remaining entries give
the labels of the rows containing dots in that column. Iterating the
operation yields all the cells in $\mathcal{T}$ containing dots. The
largest elements in cycles are column labels, the rest are row
labels.

\begin{example} \label{ex:tableau-by-columns}
Let $\pi=36187425$. Then we have
\begin{center}
\begin{tabular}{ccc}
\begin{tabular}{l|cccccccc|l}
 $i$   &1&2&3&4&5&6&7&8& cycle\\\hline
 $\pi$ &3&\textbf{6}&1&\textbf{8}&7&4&2&\textbf{5}&$(248)$\\
 $\pi'$&\textbf{3}&\textbf{5}&1&\textbf{6}&\textbf{7}&4&\textbf{2}
&8&$(12457)$\\
 $\pi''$&2&3&1&\textbf{5}&\textbf{6}&\textbf{4}&7&8&$(456)$\\
 $\pi'''$&\textbf{2}&\textbf{3}&\textbf{1}&4&5&6&7&8&$(123)$\\
 $\epsilon$   &1&2&3&4&5&6&7&8&
\end{tabular}
& \qquad &
$\mathcal{T}=
\begin{picture}(60,60)(0,54)

\setlength{\unitlength}{2pt}

\thinlines
\multiput(10,50)(10,0){4}{\line(0,-1){40}}
\put(50,50){\line(0,-1){20}}
\multiput(10,50)(0,-10){3}{\line(1,0){40}}
\multiput(10,20)(0,-10){2}{\line(1,0){30}}

\thicklines
\put(50,50){\line(0,-1){20}}
\put(50,30){\line(-1,0){10}}
\put(40,30){\line(0,-1){20}}
\put(40,10){\line(-1,0){30}}

\multiput(15,35)(0,-10){2}{\circle*{2}}
\multiput(25,45)(0,-10){4}{\circle*{2}}
\multiput(35,25)(0,-10){2}{\circle*{2}}
\multiput(45,45)(0,-10){2}{\circle*{2}}

\put(5,43){\small 1}
\put(5,33){\small 2}
\put(5,23){\small 4}
\put(5,13){\small 5}

\put(14,53){\small 8}
\put(24,53){\small 7}
\put(34,53){\small 6}
\put(44,53){\small 3}

\end{picture}
$
\end{tabular}
\end{center}
so $\pi=(123)(456)(12457)(248)$, and the path
$P(\pi)={\footnotesize
\begin{pmatrix}
1&2&3&4&5&6&7&8\\
v&v&h&v&v&h&h&h
\end{pmatrix}
}
$,
where $v$ and $h$ denote the vertical and horizontal edges,
respectively, so $\mathcal{T}=\Phi^{-1}(\pi)$ has dots in cells
labeled $(1,3)$, $(2,3)$, $(4,6)$, $(5,6)$, $(1,7)$,
$(2,7)$, $(4,7)$, $(5,7)$, $(2,8)$, $(4,8)$. Note that row labels
increase from top to bottom, while column labels increase from right 
to left.
\end{example}

\begin{remark} \label{rem:cycles-columns-induced}
Note that the column rule implies that there are no 1-cycles in 
this decomposition of $\pi$.
Likewise, the 1-hinge rule implies that if a cycle $c$ in the
product as above contains an element $i$, and another cycle $c'$ to
its left contains elements $j_1,j_2$ such that $j_1<i<j_2$, then
$c'$ also contains $i$.
\end{remark}

\begin{definition} \label{def:column-decomposition}
We call the (unique) representation of a permutation $\pi$ as a
product of increasing cycles subject to the conditions in Remark
\ref{rem:cycles-columns-induced} the \emph{column decomposition} of
$\pi$.
\end{definition}

Likewise, it is easy to see that we can determine $\mathcal{T}$ by
rows from top to bottom by decomposing $\mathcal{T}$ as a product of
\emph{decreasing} cycles. In this case, however, we may have
1-cycles, which will correspond to rows without dots. Here, at each
step, we will need to find positions of the successive
non-fixed-point \emph{right-to-left minima} of the subsequence of
$\pi$ consisting of values smaller than $\pi(s)$, where $s$ is the
position of the leftmost non-fixed point. After arriving at the
identity permutation, we add the remaining elements as fixed points.

\begin{example} \label{ex:tableau-by-rows}
Given $\pi=36187425$, we obtain $\pi=(765)(8764)(8732)(731)$ 
similarly to Example \ref{ex:tableau-by-columns}.
\end{example}

\begin{remark} \label{rem:cycles-rows-induced}
As in Remark \ref{rem:cycles-columns-induced}, note that the 1-hinge
rule implies that if an element $i$ occurs in some cycle $c$ in the
product as above, and another cycle $c'$ to the left of $c$ contains
elements $j_1,j_2$ such that $j_1>i>j_2$, then $c'$ also contains
$i$.
\end{remark}

\begin{definition} \label{def:row-decomposition}
We call the (unique) representation of a permutation $\pi$ as a
product of decreasing cycles subject to the conditions in Remark
\ref{rem:cycles-rows-induced} the \emph{row decomposition} of $\pi$.
\end{definition}

\section{Tableaux of restricted permutations}

Now we need to define the notion of a pattern. A \emph{(classical)
permutation pattern} is an order-isomorphism type of a sequence of
totally ordered letters. Such a sequence is then referred to as an
\emph{occurrence} or \emph{instance} of that pattern. For example,
$\sigma=214653$ contains instances of pattern
$\pi=2\text{-}3\text{-}1$ at subsequences $463$ and $453$ but no
instance of pattern $\tau=3\text{-}1\text{-}2$. In this case, we say
that $\sigma$ \emph{contains} $\pi$ and \emph{avoids} $\tau$. The
dashes, which are often dropped when referring to classical
patterns, are used to indicate that the terms involved in an
occurrence of the pattern may be separated by an arbitrary number of
other terms. A \emph{generalized permutation pattern} (introduced by
\cite{BaSt}) is a pattern where some letters adjacent in a pattern
must also be adjacent in a containing permutation. Such adjacent
pairs in a pattern are then \emph{not} separated by a dash. For
example, $453$ is an instance of pattern $2\text{-}31$ in $214653$,
while $463$ is an instance of pattern $23\text{-}1$, but not
$2\text{-}31$. If $\pi$ is a pattern and $\sigma$ is a permutation,
we write $(\pi)\sigma$ for the number of occurrences of $\pi$ in
$\sigma$. We also write
$\mathfrak{S}_n(\pi)=\{\sigma\in\mathfrak{S}_n\,|\,(\pi)\sigma=0\}$.
One of the earliest results concerning patterns is that
$|\mathfrak{S}_n(\pi)|=C_n$, the $n$th Catalan number, for any
classical pattern $\pi\in\mathfrak{S}_3$. The same holds for $\pi=
2\text{-}31$ and $\pi=31\text{-}2$ since
$\mathfrak{S}_n(2\text{-}31)=\mathfrak{S}_n(2\text{-}3\text{-}1)$
and
$\mathfrak{S}_n(31\text{-}2)=\mathfrak{S}_n(3\text{-}1\text{-}2)$.
See \cite{Bona-book} for more on patterns and pattern avoidance.

We will now show a nice application of row decomposition. In
\cite{StWi}, another bijection $\Psi:\mathfrak{S}_n\to
\mathfrak{S}_n$ is given (to be described later on in this paper)
that implies that the permutation tableaux containing a single 1 in
each column (i.e. the fewest possible number of 1s) are the images
of permutations avoiding pattern $2\text{-}31$, and thus are counted
by the Catalan number $C_n=\frac{1}{n+1}\binom{2n}{n}$. Open Problem
6 in \cite{StWi} asks for ``a bijection from these tableaux to any
well-known set of objects enumerated by Catalan numbers''. We will
give such a bijection using Remark \ref{rem:cycles-rows-induced}.

\begin{theorem} \label{thm:tableaux-noncrossing}
There is a natural bijection between $k$-row permutation tableaux
with $n-k$ 1s (i.e. a single 1 per column) and noncrossing
partitions of $[n]$ with $k$ blocks.
\end{theorem}

Hence, the number of permutations corresponding to these tableaux
that have $k$ weak excedances is the number of partitions of $[n]$
with $k$ blocks, i.e. the Narayana number
$N(n,k)=\frac{1}{n}\binom{n}{k}\binom{n}{k-1}$ for any $k\in[n]$.

\begin{proof}
Suppose that $\mathcal{T}$ is a tableau with a single 1 in each
column, and let $\pi=\Phi(\mathcal{T})\in \mathfrak{S}_n$. Consider
the row decomposition of $\pi$. Since each column contains a single
1, each column label occurs in only one of the cycles. Each row
label also occurs in a single cycle. Thus, every label occurs once
in the row decomposition of $\pi$, hence the cycles in the row
decomposition of $\pi$ are mutually disjoint, and the underlying
sets for these cycles form a set partition $\Pi$ of $[n]$. Suppose
that there cycles $\gamma_1\ne\gamma_2$ in the row decomposition of
$\pi$ such that $\gamma_1$ contains elements $a,c$ and $\gamma_2$
contains elements $b,d$ such that $a>b>c>d$. Thus, if $\gamma_2$ is
to the left of $\gamma_1$, then $\gamma_2$ also contains $b$, and if
$\gamma_1$ is to the left of $\gamma_2$, then $\gamma_1$ also
contains $c$. This contradicts the fact that each element in $[n]$
must occur in a single cycle of the row decomposition of $\pi$, so
$\Pi$ is noncrossing.
\end{proof}

We note that the set $E_n$ of permutations whose cycle decomposition
corresponds to noncrossing partitions of $n$ occurs in \cite{BEM} as
the set of sequences of halves of even-valued entries of
$3\text{-}1\text{-}4\text{-}2$ avoiding Dumont permutations of the
second kind. Note also that $E_n$ is exactly the set of permutations
whose column decomposition contains only 2-cycles.

The next theorem describes the tableaux of $3\text{-}2\text{-}1$
avoiding permutations.

\begin{theorem} \label{thm:321-avoiding-tableaux}
The tableaux of $3\text{-}2\text{-}1$ avoiding permutations are
exactly those whose rows and columns are all nondecreasing from left
to right and from top to bottom, respectively.
\end{theorem}

\begin{proof}
Note that $\pi$ is a $3\text{-}2\text{-}1$ avoiding permutation if
and only if each element of $\pi$ is either a left-to-right maximum
or a right-to-left-minimum, i.e. if and only if $\pi$ is the
identity or a union of two nondecreasing subsequences. Again,
without loss of generality assume that $n$ is not a fixed point of
$\pi$, i.e. $\pi(n)<n$. Suppose that the first column contains a 0
at row $l$ that is underneath a 1 at row $i$. Then as in the proof
of Theorem \ref{thm:left-column-dots}, we have
$i<l\le\pi(l)<\pi(i)$. Hence, to avoid an occurrence of
$3\text{-}2\text{-}1$ in $\pi$, we must have $\pi(m)>\pi(l)\ge l$
for all $m>l$, so $\pi$ must have at least $n-l+1$ values greater
than $\pi(l)\ge l$ ($\pi(i)$ and all $\pi(m)$ for $m>l$), which is
impossible. Therefore, the leftmost column of $\pi$ must have all 0s
atop all 1s. Moreover, if the rows containing a dot in the leftmost
column are labeled $i_j\ (1\le j\le r)$, then the sequence
$\{\pi(i_j)\}_{j=1}^{r}$, is increasing.

Let $\pi'$ be the permutation defined as in Theorem
\ref{thm:left-column-dots}. Assume that the sequence of wex tops of
$\pi'$ (see Definition \ref{def:wex}) is increasing. We have
$\pi(l)=\pi'(l)$ for a wex bottom $l<i_1$, as well as
$\pi(i_j)=\pi'(i_{j+1})$ for $1\le j<r$, and $\pi(i_r)=\pi'(n)=n$,
so the sequence of wex tops of $\pi$ is increasing if the sequence
of wex tops of $\pi'$ is increasing. Since the identity tableau has
no columns and increasing wex tops, we see by induction that if each
column of $\pi$ has all 0s atop all 1s, then the sequence of wex
tops of $\pi$ is increasing.

Thus, if $\pi$ is $3\text{-}2\text{-}1$ avoiding then the sequence
of wex tops of $\pi$ is increasing and each column of the tableau of
$\pi$ has all 0s atop of all 1s. Likewise, $\pi$ is the sequence of
nonwex bottoms of $\pi$ is increasing and each row of the tableau of
$\pi$ has all 0s to the left of all 1s. Conversely, if both wex tops
and nonwex bottoms are increasing, then $\pi$ is a union of at most
2 subsequences and hence avoids $3\text{-}2\text{-}1$.
\end{proof}

Now we will describe a bijection $\Psi$ on $\mathfrak{S}_n$ that
translates certain pattern statistics on permutations into
\emph{alignment} and \emph{crossing} statistics (see \cite{StWi}) of
the tableaux of their images. Given a permutation $\pi\in S_n$ and
we define the \emph{right embracing number} $\textsc{remb}(\pi(i))$
of the entry $\pi(i)$ as the number of instances of $2\text{-}31$ in
$\pi$ that start with $\pi(i)$. We also define the sequence
$\textsc{remb}(\pi)=\{\textsc{remb}(\pi(i))\,|\,i\in[n]\}$. If
$\pi(i)>\pi(i+1)$, then $\pi(i)\pi(i+1)$ is a \emph{descent} of
$\pi$, $\pi(i)$ a \emph{descent top} and $\pi(i+1)$ is a
\emph{descent bottom}. Similarly, if $\pi(i)<\pi(i+1)$, then
$\pi(i)$ is a \emph{non-descent bottom} and $\pi(i+1)$ is a
\emph{non-descent top}. Note that a permutation is uniquely
determined by the sets of its descents tops and descent bottoms and
the sequence of its right embracing numbers \cite{StWi}. For
$\pi\in\mathfrak{S}_n$, define the following sets:
\begin{alignat*}{2}
\textsc{db}(\pi)&=\text{set of descent bottoms of $\pi$}, &\qquad
\textsc{ndt}(\pi)&=\text{set of non-descent tops of $\pi$},\\
\textsc{dt}(\pi)&=\text{set of descent tops of $\pi$}, &\qquad
\textsc{ndb}(\pi)&=\text{set of non-descent bottoms of $\pi$},\\
\textsc{wexb}(\pi)&=\text{set of wex bottoms of $\pi$}, &\qquad
\textsc{nwext}(\pi)&=\text{set of nonwex tops of $\pi$},\\
\textsc{wext}(\pi)&=\text{set of wex tops of $\pi$}, &\qquad
\textsc{nwexb}(\pi)&=\text{set of nonwex bottoms of $\pi$}.
\end{alignat*}
Note that $\textsc{ndt}(\pi)=[n]\setminus\textsc{db}(\pi)$,
$\textsc{ndb}(\pi)=[n]\setminus\textsc{dt}(\pi)$,
$\textsc{nwext}(\pi)=[n]\setminus\textsc{wexb}(\pi)$,
$\textsc{nwexb}(\pi)=[n]\setminus\textsc{wext}(\pi)$.

The permutation $\sigma=\Psi(\pi)$ is defined as follows. We set
\begin{alignat*}{2}
\textsc{wexb}(\sigma)&=\{a+1\,|\,a\in\textsc{db}(\pi)\}\cup\{1\},
&\qquad
\textsc{nwext}(\sigma)&=\{a+1\,|\,a\in\textsc{ndt}(\pi),\ a\ne n\},\\
\textsc{wext}(\sigma)&=\{a-1\,|\,a\in\textsc{dt}(\pi)\}\cup\{n\},
&\qquad \textsc{nwexb}(\sigma)&=\{a-1\,|\,a\in\textsc{ndb}(\pi),\
a\ne 1\}
\end{alignat*}
Next we use $\textsc{remb}(\pi)$ to determine the bijections
$\textsc{wexb}(\sigma)\to\textsc{wext}(\sigma)$ and
$\textsc{nwext}(\sigma)\to\textsc{nwexb}(\sigma)$ that together form
$\Psi$. First, we find $a=\max\textsc{wexb}(\sigma)$, then find the
$(\textsc{remb}_{\pi}(a)+1)$-st smallest element in
$\textsc{wext}(\sigma)$ among those greater than or equal to $a$,
call it $b$. Then $b=\sigma(a)$. Next we delete $a$ from
$\textsc{wexb}(\sigma)$ and $b$ from $\textsc{wext}(\sigma)$ and
iterate this process until there are no more wex bottoms left.
Similarly, we find $c=\min\textsc{nwext}(\sigma)$, then find
$(\textsc{remb}_{\pi}(c)+1)$-st largest element in
$\textsc{nwexb}(\sigma)$ among those less than $c$, call it $d$.
Then $d=\sigma(c)$. Next we delete $c$ from $\textsc{nwext}(\sigma)$
and $d$ from $\textsc{nwexb}(\sigma)$ and iterate this process until
there are no more nonwex tops left.

Given a pattern $\tau$, let $(\tau)\sigma$ be the number of occurrences of $\tau$ in $\sigma$. Also, let $\mathrm{des}\,\sigma$ be the number of descents of $\sigma$, and let $\mathrm{wex}\,\pi$ be the number of weak excedances of $\pi$. It was shown in \cite{StWi} that for $\pi=\Psi(\sigma)$, we have
\[
\begin{split}
\textrm{des}\,\sigma&=\textrm{wex}\,\pi-1\\
(31\text{-}2)\sigma&=A_{EE}(\pi)+A_{NN}(\pi)\\
(21\text{-}3)\sigma+(3\text{-}21)\sigma-\binom{\textrm{des}\,\sigma}{2}&=A_{EN}(\pi)\\
(2\text{-}31)\sigma&=C_{EE}(\pi)+C_{NN}(\pi)\\
(1\text{-}32)\sigma+(32\text{-}1)\sigma-\binom{\textrm{des}\,\sigma}{2}&=A_{NE}(\pi)
\end{split}
\]

This implies (via the reverse complement map applied to $\sigma$) that the following pairs of statistics are equidistributed on permutations:
\begin{enumerate}
\item $A_{EE}+A_{NN}$ and $C_{EE}+C_{NN}$,
\item $A_{EN}$ and $A_{NE}$.
\end{enumerate}

It is surprisingly easy to show that the statistics $A_{EN}$ and $A_{NE}$ are equidistributed. Indeed, it is a straightforward exercise to show that the map $irc=i\circ r\circ c$ (inverse of reversal of complement, or reflection across the antidiagonal of the permutation diagram) preserves $\mathrm{wex}$, $A_{EE}$, $A_{NN}$, $C_{EE}$, $C_{NN}$, and exchanges $A_{EN}$ and $A_{NE}$.

In other words, $\pi'=irc(\pi)$ means that $\pi'(i)=j$ if and only if $\pi(n+1-j)=n+1-i$, so $\pi=irc(\pi')$ as well, and we have
\begin{equation} \label{eq:aen-ane-equidistributed}
\begin{split}
\mathrm{wex}\,\pi'&=\mathrm{wex}\,\pi\\
A_{EE}(\pi)&=A_{EE}(\pi')\\
A_{NN}(\pi)&=A_{NN}(\pi')\\
C_{EE}(\pi)&=C_{EE}(\pi')\\
C_{NN}(\pi)&=C_{NN}(\pi')\\
A_{EN}(\pi)&=A_{NE}(\pi')\\
A_{NE}(\pi)&=A_{EN}(\pi')
\end{split}
\end{equation}

Recall that \cite{StWi} showed that
$\Psi(\mathfrak{S}_n(2\text{-}31))$ is the set of permutations whose
tableaux have a single 1 per column. Here we establish another
pattern-related result concerning $\Psi$. Recall 
that $\mathfrak{S}_n(2\text{-}31)=\mathfrak{S}_n(2\text{-}3\text{-}
1)$ and $\mathfrak{S}_n(31\text{-}2)=\mathfrak{S}_n(3\text{-}1\text
{-}2)$.

\begin{theorem} \label{thm:psi-312}
$\Psi(\mathfrak{S}_n(31\text{-}2))=\mathfrak{S}_n(3\text{-}2\text{-}
1)$,
i.e. $\Psi(\mathfrak{S}_n(31\text{-}2))$ is the set of permutations
whose tableaux have nondecreasing rows and columns.
\end{theorem}
\begin{proof}
Let $\pi\in\mathfrak{S}_n$, and let $\sigma=\Psi(\pi)$. It is shown
in \cite{StWi} that the occurrences of $3\text{-}12$ in  $\pi$
correspond to pairs $i,j\in[n]$ such that
$j<i\le\sigma(i)<\sigma(j)$ or $\sigma(j)<\sigma(i)<i<j$. In other
words, $i$ and $j$ are both wex bottoms or both nonwex tops (i.e.
$\sigma(i)$ and $\sigma(j)$ are both wex tops or both nonwex
bottoms), and $\sigma$ has an inversion at positions $(i,j)$.
Therefore, $\pi$ avoids $31\text{-}2$ if and only if the sequence of
wex tops of $\sigma$ and the sequence of nonwex bottoms of $\sigma$
both have no inversions, i.e. are increasing. In other words,
$\sigma$ is a union of at most two increasing subsequences, that is
$\sigma\in\mathfrak{S}_n(3\text{-}2\text{-}1)$.
\end{proof}

\section{Essential 1s}

\begin{definition} \label{def:ess-1s}
Given a permutation tableau $\mathcal{T}$, we call the topmost 1 in
each column and the leftmost 1 in each row an \emph{essential} 1. If
a 1 is both the leftmost 1 in its row and the topmost 1 in its
column, then we call it \emph{doubly essential}. Let
$\mathrm{ess}(\mathcal{T})$ and $\mathrm{dess}(\mathcal{T})$ be the
number of essential 1s and doubly essential 1s in $\mathcal{T}$ and 
let $\mathrm{emptyrows}(\mathcal{T})$ be the number of rows in 
$\mathcal{T}$ that have no 1s.
\end{definition}

Note that for any tableau $\mathcal{T}$ in a $k\times(n-k)$
rectangle, we have
\[
\mathrm{ess}(\mathcal{T})+\mathrm{dess}(\mathcal{T})+\mathrm
{emptyrows}(\mathcal{T})=k+(n-k)=n
\]
by a simple sieve argument. Since each nonessential 1 is induced
by a pair of essential 1s, each tableau is determined by its
essential 1s. It was conjectured in \cite{StWi} that distribution of
the bistatistic $(\mathrm{ess},\mathrm{rows})$ of the number of
essential 1s and the number of rows on permutation tableaux is the
same as that of the bistatistic $(\mathrm{cycles},\mathrm{wex})$ of
the number of cycles and the number of weak excedances on
permutations. We will show the following.

\begin{theorem} \label{thm:ess-1s-distribution}
The bistatistic $(\mathrm{dess}+\mathrm{emptyrows},\mathrm{rows})$ 
on tableaux bounded by a path of length $n$, and the bistatistic $
(\mathrm{cycles},\mathrm{wex})$ on permutations in $\mathfrak{S}_n$ 
have the same distribution.
\end{theorem}

Since each zero row corresponds to a fixed point, it is enough to 
show that the theorem is true on derangements, so that the 
corresponding tableaux have a 1 in each row. In other words, we need 
to prove the following theorem.

\begin{theorem} \label{thm:ess-1s-dist-derangements}
The bistatistic $(\mathrm{dess},\mathrm{rows})$ on tableaux bounded 
by a path of length $n$ that have a 1 in each row, and the 
bistatistic $(\mathrm{cycles},\mathrm{wex})$ on derangements in $
\mathfrak{S}_n$ have the same distribution.
\end{theorem}

Consideration of essential 1s suggests a different type of tableau,
which we will call a \emph{bare tableau} or a \emph{0-hinge
tableau}, that results from removing all nonessential 1s from a
permutation tableau. A bare tableau is defined almost the same way
as a permutation tableau, except that the 1-hinge property is
replaced with the corresponding 0-hinge property:
\begin{description}
\item[(0-hinge)]
A cell in $Y_\lambda$ with a $1$ above it in the same column and a
$1$ to its left in the same row must contain a $0$.
\end{description}

We also define two maps $\theta$ and $\phi$ (see (\ref{eq:bare})) on tableaux of the same shape that result in removal and filling of all nonessential 1s, respectively, as well as a map $\Theta$ from bare tableaux to permutations, similar to $\Phi$, i.e. given by southeast paths from the northwestern to the southeastern boundary of the tableau that switch direction at each 1. In the diagram below, empty 
circles denote essential 1s and double circles denote doubly essential 1s.

\begin{equation} \label{eq:bare}
\begin{array}{ccc} 
\phi(\mathcal{B})=\mathcal{T}=
\begin{picture}(60,60)(0,48)
\setlength{\unitlength}{2pt}
\thinlines
\multiput(10,50)(10,0){4}{\line(0,-1){40}}
\put(50,50){\line(0,-1){20}}
\multiput(10,50)(0,-10){3}{\line(1,0){40}}
\multiput(10,20)(0,-10){2}{\line(1,0){30}}
\thicklines
\put(50,50){\line(0,-1){20}}
\put(50,30){\line(-1,0){10}}
\put(40,30){\line(0,-1){20}}
\put(40,10){\line(-1,0){30}}
\thinlines
\multiput(15,35)(10,10){2}{\circle{4}}
\multiput(15,35)(10,10){2}{\circle{2}}
\multiput(15,25)(20,0){2}{\circle{3}}
\put(25,15){\circle{3}}
\put(45,45){\circle{3}}
\multiput(25,35)(0,-10){2}{\circle*{2}}
\put(35,15){\circle*{2}}
\put(45,35){\circle*{2}}
\put(5,43){\small 1}
\put(5,33){\small 2}
\put(5,23){\small 4}
\put(5,13){\small 5}
\put(14,53){\small 8}
\put(24,53){\small 7}
\put(34,53){\small 6}
\put(44,53){\small 3}
\end{picture}
\hskip 0.5in &\longleftrightarrow &\theta(\mathcal{T})=\mathcal{B}=
\begin{picture}(60,60)(0,48)
\setlength{\unitlength}{2pt}
\thinlines
\multiput(10,50)(10,0){4}{\line(0,-1){40}}
\put(50,50){\line(0,-1){20}}
\multiput(10,50)(0,-10){3}{\line(1,0){40}}
\multiput(10,20)(0,-10){2}{\line(1,0){30}}
\thicklines
\put(50,50){\line(0,-1){20}}
\put(50,30){\line(-1,0){10}}
\put(40,30){\line(0,-1){20}}
\put(40,10){\line(-1,0){30}}
\thinlines
\multiput(15,35)(10,10){2}{\circle{4}}
\multiput(15,35)(10,10){2}{\circle{2}}
\multiput(15,25)(20,0){2}{\circle{3}}
\put(25,15){\circle{3}}
\put(45,45){\circle{3}}
\put(5,43){\small 1}
\put(5,33){\small 2}
\put(5,23){\small 4}
\put(5,13){\small 5}
\put(14,53){\small 8}
\put(24,53){\small 7}
\put(34,53){\small 6}
\put(44,53){\small 3}
\end{picture}
\\[42pt]
& &
=
\begin{picture}(30,30)(0,36)
\setlength{\unitlength}{2pt}
\thinlines
\multiput(10,30)(10,0){2}{\line(0,-1){20}}
\put(30,30){\line(0,-1){10}}
\multiput(10,30)(0,-10){2}{\line(1,0){20}}
\put(10,10){\line(1,0){10}}
\thicklines
\put(30,30){\line(0,-1){10}}
\put(30,20){\line(-1,0){10}}
\put(20,20){\line(0,-1){10}}
\put(20,10){\line(-1,0){10}}
\thinlines
\put(15,25){\circle{4}}
\put(15,25){\circle{2}}
\multiput(15,15)(10,10){2}{\circle{3}}
\put(5,23){\small 1}
\put(5,13){\small 5}
\put(14,33){\small 7}
\put(24,33){\small 3}
\end{picture}
\hskip 0.5in
\oplus
\begin{picture}(30,30)(0,36)
\setlength{\unitlength}{2pt}
\thinlines
\multiput(10,30)(10,0){3}{\line(0,-1){20}}
\multiput(10,30)(0,-10){3}{\line(1,0){20}}
\thicklines
\put(30,30){\line(0,-1){20}}
\put(30,10){\line(-1,0){20}}
\thinlines
\put(15,25){\circle{4}}
\put(15,25){\circle{2}}
\multiput(15,15)(10,0){2}{\circle{3}}
\put(5,23){\small 2}
\put(5,13){\small 4}
\put(14,33){\small 8}
\put(24,33){\small 6}
\end{picture}
\\[24pt]
\Phi(\mathcal{T})=36187425=(13)(264857)\;&\longleftrightarrow&
\;\Theta(\mathcal{B})=(1375)(2648)=36781452
\end{array}
\end{equation}

It is easy to see that the diagram $D(\mathcal{B})$ of a bare
tableau is a binary forest, since every vertex may only have a
single edge connecting it to a vertex above or to the left of it and
at most two edges to vertices south and east of it. Indeed, if there
are edges south and east to the same vertex, then $\mathcal{B}$ has
a 1-hinge at that vertex, which is impossible. Moreover, the doubly
essential 1s are exactly the roots of those binary trees.

The resulting trees are labeled as follows. The root has the two
labels of the row and column of the cell that contains it. Each
nonroot vertex has a single label: the column (resp. row) label of
the cell containing it if that vertex is reached by an edge south
(resp. east) when traveling from a root vertex. Call a nonroot
vertex a \emph{left} son if it get a row label, and a \emph{right}
son if it gets a column label. Then the following property is easy
to see: each left son has the least label in its subtree, and each 
right son has the greatest label in its subtree. Also, the root gets 
the least and the greatest label in each tree, and thus has the 
properties of both a left son and a right son.

Thus, each bare tableau $\mathcal{B}$ can be decomposed into several
tableaux each of which corresponds to a single binary tree in
$D(\mathcal{B})$ labeled as described above.
\begin{equation} \label{eq:trees}
\mathcal{B}=
\begin{picture}(60,60)(0,48)
\setlength{\unitlength}{2pt}
\thinlines
\multiput(10,50)(10,0){4}{\line(0,-1){40}}
\put(50,50){\line(0,-1){20}}
\multiput(10,50)(0,-10){3}{\line(1,0){40}}
\multiput(10,20)(0,-10){2}{\line(1,0){30}}
\thicklines
\put(50,50){\line(0,-1){20}}
\put(50,30){\line(-1,0){10}}
\put(40,30){\line(0,-1){20}}
\put(40,10){\line(-1,0){30}}
\thinlines
\multiput(15,35)(10,10){2}{\circle{4}}
\multiput(15,35)(10,10){2}{\circle{2}}
\multiput(15,25)(20,0){2}{\circle{3}}
\put(25,15){\circle{3}}
\put(45,45){\circle{3}}
\put(5,43){\small 1}
\put(5,33){\small 2}
\put(5,23){\small 4}
\put(5,13){\small 5}
\put(14,53){\small 8}
\put(24,53){\small 7}
\put(34,53){\small 6}
\put(44,53){\small 3}
\end{picture}
\hskip 0.5 in
=
\begin{picture}(30,30)(0,36)
\setlength{\unitlength}{2pt}
\thinlines
\multiput(10,30)(10,0){2}{\line(0,-1){20}}
\put(30,30){\line(0,-1){10}}
\multiput(10,30)(0,-10){2}{\line(1,0){20}}
\put(10,10){\line(1,0){10}}
\thicklines
\put(30,30){\line(0,-1){10}}
\put(30,20){\line(-1,0){10}}
\put(20,20){\line(0,-1){10}}
\put(20,10){\line(-1,0){10}}
\thinlines
\put(15,25){\circle{4}}
\put(15,25){\circle{2}}
\multiput(15,15)(10,10){2}{\circle{3}}
\put(5,23){\small 1}
\put(5,13){\small 5}
\put(14,33){\small 7}
\put(24,33){\small 3}
\end{picture}
\hskip 0.5in
\oplus
\begin{picture}(30,30)(0,36)
\setlength{\unitlength}{2pt}
\thinlines
\multiput(10,30)(10,0){3}{\line(0,-1){20}}
\multiput(10,30)(0,-10){3}{\line(1,0){20}}
\thicklines
\put(30,30){\line(0,-1){20}}
\put(30,10){\line(-1,0){20}}
\thinlines
\put(15,25){\circle{4}}
\put(15,25){\circle{2}}
\multiput(15,15)(10,0){2}{\circle{3}}
\put(5,23){\small 2}
\put(5,13){\small 4}
\put(14,33){\small 8}
\put(24,33){\small 6}
\end{picture}
\hskip 0.5 in \longleftrightarrow\;
\begin{picture}(80,10)(0,10)
 \thicklines
 \put(20,20){\circle*{3}}
 \put(10,10){\circle*{3}}
 \put(30,10){\circle*{3}}
 \put(20,20){\line(-1,-1){10}}
 \put(20,20){\line( 1,-1){10}}
 \put(70,20){\circle*{3}}
 \put(60,10){\circle*{3}}
 \put(70, 0){\circle*{3}}
 \put(70,20){\line(-1,-1){10}}
 \put(60,10){\line( 1,-1){10}}
 \put(11,15){\scriptsize $l$}
 \put(27,15){\scriptsize $r$}
 \put(61,15){\scriptsize $l$}
 \put(67,5){\scriptsize $r$}
 \put(15,25){\scriptsize 1,7}
 \put(8,0){\scriptsize 3}
 \put(28,0){\scriptsize 5}
 \put(65,25){\scriptsize 2,8}
 \put(58,0){\scriptsize 4}
 \put(68,-10){\scriptsize 6}
\end{picture}
\end{equation}

\vskip 0.35in

The permutation $\tau=\Theta(\mathcal{B})$ is obtained by traversing
each labeled binary tree as in (\ref{eq:trees}) according to the following algorithm:
\begin{enumerate}
\item Start from the smallest label (at the root) along the left edge, if
possible. If there is no left child, this is the first return to the
root (see last step).

\item At each step, start at the previous vertex and
\begin{enumerate}

\item try to move away from the root alternating unused left and right
edges as far as possible;

\item otherwise (if there are no such edges) move towards the root
along the same-side edges as far as possible.

\end{enumerate}

\item The label of the end vertex of this path is the next term in the
cycle.

\item At the first return to the root, the next term is the largest
label at the root. At the second return to the root (and when the
root has no right child), the cycle is complete.
\end{enumerate}

For example, the traversal of the trees in (\ref{eq:trees}) yields
the cycles $(1375)$ and $(2648)$ as desired.

Note that the removal of nonessential 1s leaves the labels of rows
and columns the same, and rows still correspond to weak excedances
(in particular, zero rows correspond to fixed points), while columns
correspond to deficiencies. Thus, to prove Theorem
\ref{thm:ess-1s-distribution} we first need to prove the following.

\begin{theorem} \label{thm:cyclic}
If $\mathcal{B}$ is bare tableau with no zero rows and a single
doubly essential 1, then $\tau=\Theta(B)$ is a cyclic permutation of
length greater than 1. Moreover, $\mathcal{B}$ is uniquely
recoverable from $\tau$ in the cycle notation.
\end{theorem}

\begin{proof}
The doubly essential 1 (i.e. the root of the corresponding labeled
binary tree) must have the labels $(1,n)$ where $n$ is the length of
the southeast boundary $P$ of $\mathcal{B}$. Suppose the successive
left sons away from root the are labeled $i_1,i_2,\dots,i_r$. Then
$1<i_1<i_2<\dots<i_r<n$, and each $i_{j+1}$ is the left son of $i_j$
(letting $i_0:=1$ and $i_{r+1}:=n$). Let $D(\mathcal{B}_{i_j})$ be
the subtree of $D(\mathcal{B})$ with the right son of $i_j$ as the
root. Then, denoting the traversal of $D(\mathcal{B})$ by $\mathrm
{tr}(D(\mathcal{B}))$, we have 
\[
\mathrm{tr}(D(\mathcal{B}))=1,
\mathrm{tr}(D(\mathcal{B}_{i_1})),i_1,
\mathrm{tr}(D(\mathcal{B}_{i_2})),i_2,
\dots,\mathrm{tr}(D(\mathcal{B}_{i_r})),i_r,n,
\mathrm{tr}(D(\mathcal{B}_{n}))
\]
Each $\mathcal{B}_{i_j}$ is strictly smaller than $\mathcal{B}$, so
by inductive assumption the traversal of $\mathcal{B}_{i_j}$
contains every vertex label of $\mathcal{B}_{i_j}$ once. Hence, the
traversal of $\mathcal{B}$ contains every label of $\mathcal{B}$
once, i.e. $\tau=\Theta(\mathcal{B})$ is a cyclic permutation.

Note also that the first part of Theorem \ref{thm:left-column-dots}
(i.e. that $\pi=\pi'\circ(i_1\,i_2\,\dots\,i_r\,n)$) remains true
for bare tableaux after changing $\Phi$ to $\Theta$. Thus, to show
that $\Theta$ is a bijection we only need to prove that we can
uniquely recover the leftmost column of the bare tableau
$\mathcal{B}$ with no fixed points and a single doubly essential 1
from the cyclic permutation $\tau=\Theta(\mathcal{B})$.

It is not hard to see that because $i_j$ is a left son, every label
in $D(\mathcal{B}_{i_j})$ is larger than $i_j$. Hence, if $\rho$ is
obtained by removing the parentheses from $\tau$ in the
cycle notation starting with the least element (e.g. 1), then $\{i_r, i_{r-1},\dots,i_1,1\}$ are the right-to-left minima of the prefix of $\rho$ before $n$.
\end{proof}

Now Theorem \ref{thm:cyclic} clearly implies Theorem \ref
{thm:ess-1s-dist-derangements}, and hence, Theorem
\ref{thm:ess-1s-distribution}. Note that we can similarly prove that
the column labels of vertices in the top row (i.e. row 1) of
$\mathcal{B}$ are the right-to-left maxima of the suffix of $\rho$
starting with the greatest element (e.g. $n$). 

\begin{example} \label{ex:bare-tableau-from-permutation}
If $\tau=(2648)$, then $\rho=2648$, so $\mathcal{B}=\theta(\tau)$ has top row label 2 and leftmost column label 8, so the leftmost column of 
$\mathcal{B}$ has 1s in rows labelled 4 and 2, while the top row of 
$\mathcal{B}$ only has a 1 in column labeled 8.
\end{example}

\bigskip

In conclusion, we mention some open problems related to permutation tableaux.
\begin{enumerate}
\item Describe the map $\Phi^{-1}\circ irc\circ\Phi$ (see \eqref{eq:aen-ane-equidistributed}) directly on tableaux, without mapping to permutations and then back to tableaux.

\item Define the map $B$ as ``$\Phi$ on bare tableaux'' (in other words, a permutation is obtained from a bare tableau by similarly constructed southeast paths). Construct a direct bijection between $B^{-1}(\pi)$ and $\Phi^{-1}(\pi)$.

\item Let $a_{ij}$ be the number of derangements $\pi$ such that
$B^{-1}(\pi)$ has $i$ essential 1s, and $\Phi^{-1}(\pi)$ has
$j$ essential 1s (so $n\le i,j\le 2n-1$). Let $A=[a_{ij}]_{n\times n}$. What can be said about various $a_{ij}$'s or the whole matrix $A$?
\end{enumerate}

\end{document}